\documentclass[11.05pt,a4paper]{amsart}
\usepackage{amssymb,amsmath}
\usepackage{color}
\usepackage{latexsym}
\usepackage{amsthm,amsfonts,amssymb,mathrsfs}
\usepackage{rotating}
\usepackage[leqno]{amsmath}
\usepackage{xspace}
\usepackage[all]{xy}
\usepackage{longtable}

\headsep=1cm \numberwithin{equation}{section}
\newtheorem{theorem}{Theorem}[section]
\newtheorem{lemma}[theorem]{Lemma}
\newtheorem{proposition}[theorem]{Proposition}
\newtheorem{corollary}[theorem]{Corollary}

\theoremstyle{definition}

\theoremstyle{remark}

\newtheorem{example}[theorem]{Example}
\newtheorem{question}[theorem]{Question}

\newcommand{\im}{\operatorname{im}}

\newcommand{\Assh}{\operatorname{Assh}}
\newcommand{\Spec}{\operatorname{Spec}}

\newcommand{\cd}{\operatorname{cd}}

\newcommand{\Ht}{\operatorname{ht}}

\newcommand{\id}{\operatorname{id}}

\newcommand{\Gdim}{\operatorname{G--dim}}

\newcommand{\V}{\operatorname{V}}

\newcommand{\Ext}{\operatorname{Ext}}
\newcommand{\Supp}{\operatorname{Supp}}

\newcommand{\Hom}{\operatorname{Hom}}

\newcommand{\Ann}{\operatorname{Ann}}

\newcommand{\depth}{\operatorname{depth}}

\newcommand{\lo}{\longrightarrow}
\newcommand{\fm}{\frak{m}}
\newcommand{\fp}{\frak{p}}

\newcommand{\fa}{\frak{a}}

\newenvironment{prf}[1][Proof]{\begin{proof}[\bf #1]}{\end{proof}}

\begin{document}

\author[F. Mohammadi Aghjeh Mashhad and K. Divaani-Aazar]{Fatemeh Mohammadi Aghjeh
Mashhad and Kamran Divaani-Aazar}
\title[Foxby equivalence, local duality and ...]
{Foxby equivalence, local duality and Gorenstein homological dimensions}

\address{F. Mohammadi, Science and Research Branch, Islamic Azad
University, Tehran, Iran.} \email{mohammadi\_fh@yahoo.com}

\address{K. Divaani-Aazar, Department of Mathematics, Az-Zahra University,
Vanak, Post Code 19834, Tehran, Iran-and-School of Mathematics,
Institute for Research in Fundamental Sciences (IPM), P.O. Box
19395-5746, Tehran, Iran} \email{kdivaani@ipm.ir}


\keywords {Auslander categories,  Foxby equivalence, G-dimension, generalized
local cohomology modules, Gorenstein injective dimension,  local duality, normalized
dualizing complexes, totally reflexive modules. \\
The second author was supported by a grant from IPM (No. 89130118).}

\begin{abstract} Let $(R,\fm)$ be a local ring and $(-)^{\vee}$ denote the Matlis
duality functor. We investigate the relationship between Foxby equivalence and local
duality through generalized local cohomology modules. Assume that $R$ possesses a
normalized dualizing complex $D$ and $X$ and $Y$ are two homologically bounded
complexes of $R$-modules with finitely generated homology modules. We present several
duality results for $\fm$-section complex ${\bf R}\Gamma_{\fm}({\bf R}\Hom_R(X,Y))$.
In particular, if G-dimension of $X$ and injective dimension of $Y$ are finite, then we
show that $${\bf R}\Gamma_{\fm}({\bf R}\Hom_R(X,Y))\simeq ({\bf R}\Hom_R(Y,D\otimes_
R^{{\bf L}}X))^{\vee}.$$ We deduce several applications of these duality results.
In particular, we establish Grothendieck's non-vanishing Theorem in the context of
generalized local cohomology modules.
\end{abstract}

\maketitle

\section{Introduction}

Let $(R,\fm)$ be a Cohen-Macaulay local ring with a canonical module $\omega$ and $\mathcal{P}$
(resp. $\mathcal{I}$) denote the full subcategory of finitely generated $R$-modules of finite
projective (resp. injective) dimension. By virtue of \cite[Theorem 2.9]{Sh}, there is the following
equivalence of categories
\begin{displaymath}
\xymatrix{\mathcal{P} \ar@<0.7ex>[rrr]^-{\omega\otimes_R-} &
{} & {} & \mathcal{I} \ar@<0.7ex>[lll]^-{\Hom_R(\omega,-)}}.
\end{displaymath}
Let $M$ and $N$ be two finitely generated $R$-modules and $i$ a non-negative integer. Denote the
Matlis duality functor $\Hom_R(-,E(R/\fm))$ by $(-)^{\vee}$. If $M$ has finite projective dimension,
then by Suzuki's Duality Theorem \cite[Theorem 3.5]{Su}, there is a natural isomorphism
$H_{\fm}^i(M,N)\cong \Ext_R^{\dim R-i}(N,\omega\otimes_RM)^{\vee}.$ Also, if $N$ has finite injective
dimension, then the Herzog-Zamani Duality Theorem \cite[Theorem 2.1 b)]{HZ} asserts that
$H_{\fm}^i(M,N)\cong \Ext_R^{\dim R-i}(\Hom_R(\omega,N),M)^{\vee}$. (These results can be considered
as variants of the Local Duality Theorem \cite[11.2.8]{BS} in the context of generalized local
cohomology modules.) Hence the equivalence between two subcategories $\mathcal{P}$ and $\mathcal{I}$
can be connected to local duality through generalized local cohomology modules.

Now, assume that $(R,\fm)$ is a local ring with a normalized dualizing complex $D$. By
\cite[Theorem 4.1 and Prposition 3.8 b)]{CFrH} (resp. \cite[Theorem 4.4]{CFrH}) Auslander
category $\mathcal{A}^f(R)$ (resp. $\mathcal{B}^f(R)$) consists exactly of all homologically
bounded complexes of $R$-modules with finitely generated homology modules of finite G-dimension
(resp. Gorenstein injective dimension). By Foxby equivalence, there is an equivalence of categories
\begin{displaymath}
\xymatrix{\mathcal{A}^f(R) \ar@<0.7ex>[rrr]^-{D\otimes_R^{\bf L}-}
& {} & {} & \mathcal{B}^f(R) \ar@<0.7ex>[lll]^-{{\bf R}\Hom_R(D,-)}};
\end{displaymath}
see e.g. \cite[Theorem 3.3.2 a),  b), e) and f)]{C}. Foxby equivalence between the two categories
$\mathcal{A}^f(R)$ and $\mathcal{B}^f(R)$ is a natural generalization of the above mentioned equivalence
between the subcategories $\mathcal{P}$ and $\mathcal{I}$. In view of what we saw in the first paragraph,
it is natural to ask whether Foxby equivalence can also be connected to local duality through generalized
local cohomology modules. Assume that $X$ and $Y$ are two homologically bounded complexes of $R$-modules
with finitely generated homology modules and let $i$ be an integer. The following natural questions arise:

\begin{question} Suppose that G-dimension of $X$ is finite. Is $H_{\fm}^i(X,Y) \cong \Ext_R^{-i}(Y,
D\otimes_R^{\bf L}X)^{\vee}$?
\end{question}

\begin{question} Suppose that Gorenstein injective dimension of $Y$ is finite. Is
$$H_{\fm}^i(X,Y)\cong \Ext_R^{-i}({\bf R}\Hom_R(D,Y),X)^{\vee}?$$
\end{question}

Our main aim in this paper is to answer these questions. Example 3.6 below shows that the answers of
these questions are negative in general, but by adding some extra assumptions on the complexes $X$ and
$Y,$  we can deduce our desired natural isomorphisms. Consider the following assumptions:
\begin{enumerate}
\item[a)] Projective dimension of $X$ is finite.
\item[b)] Projective dimension of $Y$ is finite.
\item[b')] Both G-dimension of $X$ and projective
dimension of $Y$ are finite.
\item[c)] Both G-dimension of $X$ and injective
dimension of $Y$ are finite.
\item[d)] Injective dimension of $Y$ is finite.
\item[e)] Injective dimension of $X$ is finite.
\item[e')] Both Gorenstein injective dimension of $Y$ and injective
dimension of $X$ are finite.
\item[f)] Both Gorenstein injective dimension of $Y$ and projective
dimension of $X$ are finite.
\end{enumerate}
We  show that each of a), b) and c) implies the natural isomorphism $${\bf R}\Gamma_{\fm}({\bf R}
\Hom_R(X,Y))\simeq ({\bf R}\Hom_R(Y,D\otimes_R^{{\bf L}} X))^{\vee},$$ and  each of d), e) and f)
implies the natural isomorphism $${\bf R}\Gamma_{\fm}({\bf R}\Hom_R(X,Y))\simeq ({\bf R}\Hom_R({\bf
R}\Hom_R(D,Y),X))^{\vee}.$$
These immediately yield our desired isomorphisms $H_{\fm}^i(X,Y)\cong \Ext_R^{-i}(Y,D\otimes_R^{\bf L}X)
^{\vee}$ and $H_{\fm}^i(X,Y)\cong \Ext_R^{-i}({\bf R}\Hom_R(D,Y),X)^{\vee}$, respectively. These duality
results are far reaching generalizations of Suzuki's Duality Theorem and the Herzog-Zamani Duality
Theorem.

We present some applications of the above duality results. First of all, we improve the main
results of \cite{HZ}; see Propositions 4.3 and 4.4 below. Then we establish an analogue of
Grothendieck's non-vanishing Theorem in the context of generalized local cohomology modules.
Let $(R,\fm)$ be a local ring and $M,N$ two finitely generated $R$-modules such that $X:=M,
Y:=N$ satisfy one of the above assumptions a),  b'), d), and e'). When $R$ is Cohen-Macaulay,
we show that $$\cd_{\fm}(M,N)=\dim R-\depth(\Ann_RN,M).$$ Finally, we give a partial generalization  
of the Intersection inequality; see Proposition 4.7 below.

\section{Prerequisites}

Throughout this paper, $R$ is a commutative Noetherian ring with nonzero identity. The
$\fm$-adic completion of an $R$-module $M$ over a local ring $(R,\fm)$ will be denoted by
$\widehat{M}$.

{\bf (2.1) Hyperhomology.} We will work within $\mathcal{D}(R)$, the derived category of
$R$-modules. The objects in  $\mathcal{D}(R)$ are complexes of $R$-modules and symbol
$\simeq$ denotes isomorphisms in this category. For a complex $$X=\cdots \lo X_{n+1}\overset
{{\partial}_{n+1}^X}\lo X_n\overset{{\partial}_n^X}\lo X_{n-1}\lo \cdots$$ in $\mathcal{D}(R)$,
its supremum and infimum are defined, respectively, by $\sup X:=\sup \{i\in \mathbb{Z}|H_i(X)\neq 0\}$
and $\inf X:=\inf \{i\in \mathbb{Z}|H_i(X)\neq 0\}$, with the usual convention that $\sup
\emptyset=-\infty$ and $\inf \emptyset=\infty$. For an integer $\ell$, $\Sigma^{\ell}X$ is the
complex $X$ shifted $\ell$ degrees to the left. Modules will be considered as complexes concentrated
in degree zero and we denote the full subcategory of complexes with homology concentrated in degree
zero by $\mathcal{D}_0(R)$. The full subcategory of complexes homologically bounded to the right
(resp. left) is denoted by $\mathcal{D}_{\sqsupset}(R)$ (resp. $\mathcal{D}_{\sqsubset}(R)$).
Also, the full subcategories of homologically bounded complexes and of complexes with finitely
generated homology modules will be denoted by $\mathcal{D}_{\Box}(R)$ and $\mathcal{D}^f(R)$,
respectively. Throughout for any two properties $\sharp$ and $\natural$ of complexes, we set
$\mathcal{D}^{\natural}_{\sharp}(R):=\mathcal{D}_{\sharp}(R)\cap \mathcal{D}^{\natural}(R)$. So
for instance, $\mathcal{D}_{\Box}^f(R)$ stands for the full subcategory of homologically bounded
complexes with finitely generated homology modules.

For any complex $X$ in $\mathcal{D}_{\sqsupset}(R)$ (resp.$\mathcal{D}_{\sqsubset}(R)$), there
is a bounded to the right (resp. left) complex $P$ (resp. $I$) consisting of projective (resp.
injective) $R$-modules which is isomorphic to $X$ in $\mathcal{D}(R)$. A such complex $P$ (resp.
$I$) is called a projective (resp. injective) resolution of $X$. A complex $X$ is said to have
finite projective (resp. injective) dimension, if $X$ possesses a bounded projective (resp.
injective) resolution. Similarly, a complex $X$ is said  to have finite flat dimension if
it is isomorphic (in $\mathcal{D}(R)$) to a  bounded complex of flat $R$-modules. The left derived
tensor product functor $-\otimes_R^{{\bf L}}\sim$ is computed by taking a projective resolution
of the first argument or of the second one. The right derived homomorphism functor
${\bf R}\Hom_R(-,\sim)$ is computed by taking a projective resolution of the first argument or
by taking an injective resolution of the second one. For any two complexes $X$ and $Y$ and any
integer $i$, set $\Ext_R^i(X,Y):=H_{-i}({\bf R}\Hom_R(X,Y))$. Let $X$ be a complex and $\fa$ an
ideal of $R$. Recall that $\Supp_RX:=\cup_{l\in \mathbb{Z}}\Supp_RH_l(X)$,  $\depth(\fa,X):=
-\sup {\bf R}\Hom_R(R/\fa,X)$ and $\dim_RX:=\sup \{\dim R/\fp-\inf X_{\fp}|\fp\in \Spec R\}$.
For any complexes $X,Y$ and $Z$, there are the following natural isomorphisms in $\mathcal{D}(R)$.
\begin{enumerate}
\item[]  {\bf Shifts:} Let $i,j$ be two integers. Then $\Sigma^{i}
X\otimes_R^{{\bf L}}\Sigma^{j}Y \simeq \Sigma^{j+i}(X\otimes_R^{{\bf
L}}Y)$ and ${\bf R}\Hom_R(\Sigma^{i}X,\Sigma^{j}Y)\simeq
\Sigma^{j-i}{\bf R}\Hom_R(X,Y)$.
\item[]  {\bf Commutativity:} $X\otimes_R^{{\bf L}}Y\simeq Y\otimes_R^{{\bf L}}X.$
\item[] {\bf Adjointness:} Let $S$ be an $R$-algebra. If $X\in \mathcal{D}_{\sqsupset}(S)$,
$Y\in \mathcal{D}(S)$ and $Z\in \mathcal{D}_{\sqsubset}(R)$, then
$${\bf R}\Hom_R(X\otimes_S^{{\bf L}}Y,Z)\simeq {\bf R}\Hom_S(X,{\bf R}\Hom_R(Y,Z)).$$
\item[] {\bf Tensor evaluation:} Assume that $X\in \mathcal{D}^f_{\sqsupset}(R)$, $Y\in
\mathcal{D}_{\Box}(R)$ and $Z\in \mathcal{D}_{\sqsupset}(R)$. If
either projective dimension of $X$ or flat dimension of $Z$ is
finite, then
$${\bf R}\Hom_R(X,Y)\otimes_R^{{\bf L}}Z\simeq {\bf R}\Hom_R(X,Y\otimes_R^{{\bf L}}Z).$$
\item[] {\bf Hom evaluation:} Assume that $X\in \mathcal{D}^f_{\sqsupset}(R)$, $Y\in
\mathcal{D}_{\Box}(R)$ and $Z\in \mathcal{D}_{\sqsubset}(R)$. If
either projective dimension of $X$ or injective dimension of $Z$ is
finite, then
$$X\otimes_R^{{\bf L}}{\bf R}\Hom_R(Y,Z)\simeq {\bf R}\Hom_R({\bf R}\Hom_R(X,Y),Z).$$
\end{enumerate}

{\bf (2.2) Gorenstein homological dimensions.}  An $R$-module $M$ is said to be \emph{
totally reflexive} if there exists an exact complex $P$ of finitely generated projective
$R$-modules such that $M\cong \im(P_0\lo P_{-1})$ and $\Hom_R(P,R)$ is exact. Also, an
$R$-module $N$ is said to be \emph{Gorenstein injective} if there exists an exact complex
$I$ of injective $R$-modules such that $N\cong \im(I_1\lo I_0)$ and $\Hom_R(E,I)$ is exact
for all injective $R$-modules $E$; see \cite{EJ}. Obviously, any finitely generated
projective $R$-module is totally reflexive and any injective $R$-module is Gorenstein
injective. A complex $X\in \mathcal{D}^f_{\Box}(R)$ is said to have finite G-dimension
if it is isomorphic (in $\mathcal{D}(R)$) to a bounded complex of totally reflexive
$R$-modules. Also, a complex $X\in \mathcal{D}_{\Box}(R)$ is said to have finite Gorenstein
injective dimension if it is isomorphic (in $\mathcal{D}(R)$) to a bounded complex of
Gorenstein injective $R$-modules.

{\bf (2.3) Auslander categories.} Let $(R,\fm)$ be a local ring. A \emph{normalized dualizing
complex} for $R$ is a complex $D\in \mathcal{D}_{\Box}^f(R)$ such that the homothety morphism
$R\longrightarrow {\bf R}\Hom_R(D,D)$ is an isomorphism in $\mathcal{D}(R)$, $D$ has finite
injective dimension and $\sup D=\dim R$. Assume that $R$ possesses a normalized dualizing complex
$D$. The Auslander category $\mathcal{A}^f(R)$ (with respect to $D$) is the full subcategory of
$\mathcal{D}^f_{\Box}(R)$ whose objects are exactly those complexes $X\in \mathcal{D}^f_{\Box}(R)$
for which $D\otimes_{R}^{{\bf L}}X\in \mathcal{D}^f_{\Box}(R)$ and the natural morphism
$\eta_X:X\lo {\bf R}\Hom_R(D,D\otimes_{R}^{{\bf L}}X)$ is an isomorphism in $\mathcal{D}(R)$.
Also, the Auslander category $\mathcal{B}^f(R)$ (with respect to $D$) is the full subcategory of
$\mathcal{D}^f_{\Box}(R)$ whose objects are exactly those complexes $X\in \mathcal{D}^f_{\Box}(R)$
for which ${\bf R}\Hom_R(D,X)\in \mathcal{D}^f_{\Box}(R)$ and the natural morphism $\varepsilon_X:
D\otimes_{R}^{{\bf L}}{\bf R}\Hom_R(D,X)\lo X$ is an isomorphism in $\mathcal{D}(R)$.
By \cite[Theorem 4.1 and Prposition 3.8 b)]{CFrH}, $\mathcal{A}^f(R)$ precisely consists of all
complexes $X\in \mathcal{D}^f_{\Box}(R)$ whose G-dimensions are finite. Also, \cite[Theorem 4.4]{CFrH}
yields that $\mathcal{B}^f(R)$ consists of all complexes $X\in \mathcal{D}^f_{\Box}(R)$ whose
Gorenstein injective dimensions are finite.

{\bf (2.4) Local cohomology.} Let $\fa$ be an ideal of $R$. The right derived functor of $\fa$-section
functor $\Gamma_{\fa}(-)={\varinjlim}_n\Hom_R(R/\fa^n,-)$ is denoted by ${\bf R}\Gamma_{\fa}(-)$.
For any complex $X\in \mathcal{D}_{\sqsubset}(R)$, the complex ${\bf R}\Gamma_{\fa}(X)\in
\mathcal{D}_{\sqsubset}(R)$ is defined by ${\bf R}\Gamma_{\fa}(X):=\Gamma_{\fa}(I)$, where $I$ is
an (every) injective resolution of $X$. Also, for any two complexes $X\in \mathcal{D}_{\sqsupset}(R)$
and $Y\in\mathcal{D}_{\sqsubset}(R)$, the generalized $\fa$-section complex ${\bf R}\Gamma_{\fa}(X,Y)$
is defined by ${\bf R}\Gamma_{\fa}(X,Y):={\bf R}\Gamma_{\fa}({\bf R}\Hom_R(X,Y))$; see \cite{Y}.
For any integer $i$, set $H^i_{\fa}(X,Y):=H_{-i}({\bf R}\Gamma_{\fa}(X,Y))$ and denote $\sup\{i\in \mathbb{Z}|H^i_{\fa}(X,Y)\neq 0\}$ by $\cd_{\fa}(X,Y)$. Let $M$ and $N$ be two $R$-modules. The notion
of generalized local cohomology modules $H^i_{\fa}(M,N):={\varinjlim}_n\Ext^{i}_{R}(M/\fa^{n}M,N)$
was introduced by Herzog in his Habilitationsschrift \cite{He}. When $M$ is finitely generated,
\cite[Theorem 3.4]{Y} yields that $H^i_{\fa}(M,N)\cong H_{-i}({\bf R}\Gamma_{\fa}(M,N))$ for all
integers $i$.

Let $\Check{C}(\underline{\fa})$ denote the $\Check{C}$ech complex on a set $\underline{\fa}=
\{x_1,x_2,\dots ,x_n\}$ of generators of $\fa$. So, by the definition, $\Check{C}(\underline{\fa})=\Check{C}(x_1)\otimes_R\cdots \otimes_R\Check{C}(x_n)$, where for each $i$, $\Check{C}(x_i)$ is the complex $0\lo R\lo R_{x_i}\lo 0$ concentrated in degrees $0$ and $-1$ in
which homomorphisms are the natural ones. For any complex $X\in \mathcal{D}_{\sqsubset}(R)$,
\cite[Theorem 1.1 iv)]{Sc} implies that ${\bf R}\Gamma_{\fa}(X)\simeq X\otimes_R^{{\bf L}}\Check{C}
(\underline{\fa}).$

\section{Duality Results}

We start by proving two lemmas which are needed in the proof of the main result of this section.

\begin{lemma} Let $(R,\fm)$ be a local ring possessing a normalized dualizing complex $D$ and $X,Y\in \mathcal{D}_{\Box}^f(R)$.
\begin{enumerate}
\item[i)] Assume that one of the following conditions holds:\\
a) either projective dimension $X$ or $Y$ is finite,\\
b) both G-dimension of $X$ and injective dimension of $Y$ are finite.\\
Then $$X\otimes_R^{{\bf L}}{\bf R}\Hom_R(Y,D)\simeq {\bf R}\Hom_R(Y,D\otimes_R^{{\bf L}}X).$$
\item[ii)] Assume that one of the following conditions holds:\\
a) either injective dimension $Y$ or $X$ is finite,\\
b) both Gorenstein injective dimension of $Y$ and projective
dimension of $X$ are finite.\\
Then  $$X\otimes_R^{{\bf L}}{\bf R}\Hom_R(Y,D)\simeq {\bf R}\Hom_R({\bf R}\Hom_R(D,Y),X).$$
\end{enumerate}
\end{lemma}

\begin{prf}  i) The case a) follows immediately by using commutativity of the bivariant functor
$-\otimes_R^{{\bf L}}\sim$ and tensor evaluation. Assume that b) holds.  Since $Y\in \mathcal{B}^f(R)$,
we have ${\bf R}\Hom_R(D,Y)\in \mathcal{D}_{\Box}^f(R)$. As $Y$ has finite injective dimension,
\cite[Theorem 3.3.2 d)]{C} and \cite[Theorem A.5.7.2]{C} imply that ${\bf R}\Hom_R(D,Y)$ has finite
projective dimension. Next, as $Y\in \mathcal{B}^f(R)$ and ${\bf R}\Hom_R(D,D)\simeq R$,
\cite[Lemma 3.3.3 b)]{C} yields that ${\bf R}\Hom_R(Y,D)\simeq {\bf R}\Hom_R({\bf R}\Hom_R(D,Y),R)$.
Now, tensor evaluation, the fact $X\in \mathcal{A}^f(R)$ and \cite[Lemma 3.3.3 b)]{C} yield that:
$$\begin{array}{llll} X\otimes_R^{{\bf L}}{\bf R}\Hom_R(Y,D)&\simeq
X\otimes_R^{{\bf L}}{\bf R}
\Hom_R({\bf R}\Hom_R(D,Y),R) \\
&\simeq {\bf R}\Hom_R({\bf R}\Hom_R(D,Y),R)\otimes_R^{{\bf L}}X \\
&\simeq {\bf R}\Hom_R({\bf R}\Hom_R(D,Y),X)\\
&\simeq {\bf R}\Hom_R({\bf R}\Hom_R(D,Y),{\bf R}\Hom_R(D,
D\otimes_R^{{\bf L}}X))\\
&\simeq {\bf R}\Hom_R(Y,D\otimes_R^{{\bf L}}X).
\end{array}
$$

ii) Assume that Gorenstein injective dimension of $Y$ is finite. Then $Y\in \mathcal{B}^f(R)$,  and so
${\bf R}\Hom_R(D,Y)\in \mathcal{D}_{\Box}^f(R)$. As, we saw in i), we have
$$X\otimes_R^{{\bf L}}{\bf R}\Hom_R(Y,D)\simeq {\bf R}\Hom_R({\bf R}\Hom_R(D,Y),R)\otimes_R^{{\bf L}}X.$$
As we mentioned above, the finiteness of injective dimension of $Y$, implies that ${\bf R}\Hom_R(D,Y)$
has finite projective dimension. Thus  b) and the first case of a) follow by tensor evaluation. It remains
to consider the second case of a). So, assume that $X$ has finite injective dimension. Now, as
$X\in \mathcal{B}^f(R)$, by using tensor evaluation and Hom evaluation, we can deduce that:
$$
\begin{array}{llll} X\otimes_R^{{\bf L}}{\bf R}\Hom_R(Y,D)&\simeq
\big(D\otimes_R^{{\bf L}}{\bf R}
\Hom_R(D,X)\big)\otimes_R^{{\bf L}}{\bf R}\Hom_R(Y,D) \\
&\simeq D\otimes_R^{{\bf L}}\big({\bf R}\Hom_R(Y,D)\otimes_R^{{\bf
L}}
{\bf R}\Hom_R(D,X)\big) \\
&\simeq D\otimes_R^{{\bf L}}{\bf R}\Hom_R\big(Y,D\otimes_R^{{\bf L}}
{\bf R}\Hom_R(D,X)\big) \\
&\simeq D\otimes_R^{{\bf L}}{\bf R}\Hom_R(Y,X) \\
&\simeq {\bf R}\Hom_R({\bf R}\Hom_R(D,Y),X).
\end{array}
$$
\end{prf}

The first assertion of the next result was already proved by Foxby \cite[Proposition 6.1]{Fo2}. For
completeness' sake, we include an easy proof for it. Recall that for a complex
$Y\in \mathcal{D}_{\sqsubset}(R)$, its injective dimension, $\id_RY$, is defined by
$$\id_RY:=\inf\{\sup\{l\in \mathbb{Z}|I_{-l}\neq 0\}|I \   \  \text{is an injective resolution of Y}\}.$$

\begin{lemma} Let $\fa$ be an ideal of $R$.  Let $X\in \mathcal{D}^f_{\sqsupset}(R)$ and  $Y\in
\mathcal{D}_{\Box}(R)$. Then ${\bf R}\Gamma_{\fa}(X,Y)\simeq {\bf R}\Hom_R(X,{\bf R}\Gamma_{\fa}(Y))$.
In particular, if $X$ is homologically bounded and not homologically trivial, then
$\cd_{\fa}(X,Y)\leq \id_RY+\sup X$.
\end{lemma}

\begin{prf} Let $\underline{\fa}$ be a generating set of $\fa$. As $\Check{C}(\underline{\fa})$ is a
bounded complex of flat $R$-modules, tensor evaluation property yields that $${\bf R}\Gamma_{\fa}(X,Y)
\simeq {\bf R}\Hom_R(X,Y)\otimes_R^{{\bf L}} \Check{C}(\underline{\fa})\simeq {\bf R}\Hom_R(X,{\bf
R}\Gamma_{\fa}(Y)).$$

Now, assume that $X$ is homologically bounded and not homologically trivial. Since $\sup X$ is an
integer, we may and do assume that $\id_RY<\infty$. So, there is a bounded complex $I$ consisting of
injective modules such that it is isomorphic to $Y$ in $\mathcal{D}(R)$ and $I_j=0$ for all $j<-\id_RY$.
One has ${\bf R}\Gamma_{\fa}(X,Y)\simeq {\bf R}\Hom_R(X,\Gamma_{\fa}(I)).$ The complex $\Gamma_{\fa}(I)$
is a bounded complex consisting of injective modules. Now by \cite[Corollary A.5.2]{C}, we have
$$\id_RY\geq \id_R\Gamma_{\fa}(I)\geq -\sup X-\inf {\bf R}\Hom_R(X,\Gamma_{\fa}(I)).$$
Thus $-\inf {\bf R}\Gamma_{\fa}(X,Y)\leq \id_RY+\sup X$, as claimed.
\end{prf}

Now, we are ready to prove the main result of this section.

\begin{theorem} Let $(R,\fm)$ be a local ring possessing a normalized dualizing complex $D$ and $X,Y\in \mathcal{D}_{\Box}^f(R)$.
\begin{enumerate}
\item[i)] Assume that one of the following conditions holds:\\
a) either projective dimension $X$ or $Y$ is finite,\\
b) both G-dimension of $X$ and injective dimension of $Y$ are finite.\\
Then $${\bf R}\Gamma_{\fm}(X,Y)
\simeq \Hom_R({\bf R}\Hom_R(Y,D\otimes_R^{{\bf L}}X),E(R/\fm)).$$
\item[ii)] Assume that one of the following conditions holds:\\
a) either injective dimension $Y$ or $X$ is finite,\\
b) both Gorenstein injective dimension of $Y$ and projective dimension of $X$ are finite.\\
Then  $${\bf R}\Gamma_{\fm}(X,Y)\simeq \Hom_R({\bf R}\Hom_R({\bf R}\Hom_R(D,Y),X),E(R/\fm)).$$
\end{enumerate}
\end{theorem}

\begin{prf} Denote the Matlis duality functor $\Hom_R(-,E(R/\fm))$ by $(-)^{\vee}$. By Local Duality
Theorem for any complex $Z\in \mathcal{D}_{\Box}^f(R)$, we know that ${\bf R}\Gamma_{\fm}(Z)\simeq ({\bf R}\Hom_R(Z,D))^{\vee}$, see e.g. \cite[Chapter V, Theorem 6.2]{Ha}. Using Lemma 3.2 and adjointness yields
that:
$$\begin{array}{llll} {\bf R}\Gamma_{\fm}(X,Y)&\simeq {\bf R}\Hom_R(X,{\bf R}\Gamma_{\fm}(Y))\\
&\simeq {\bf R}\Hom_R(X,{\bf R}\Hom_R(Y,D)^{\vee})\\
&\simeq (X\otimes_R^{{\bf L}}{\bf R}\Hom_R(Y,D))^{\vee}.
\end{array}
$$
Hence Lemma 3.1 completes the proof.
\end{prf}

Theorem 3.3 has the following immediate corollary.

\begin{corollary}
\begin{enumerate}
\item[i)] Let the situation be as in Theorem 3.3 i). Then
$$-\inf {\bf R}\Gamma_{\fm}(X,Y)=\sup {\bf R}\Hom_R(Y,D\otimes_R^{{\bf L}}X)$$ and $$-\sup {\bf
R}\Gamma_{\fm}(X,Y)=\inf {\bf R}\Hom_R(Y,D\otimes_R^{{\bf L}}X).$$
\item[ii)] Let the situation be as in Theorem 3.3 ii). Then
$$-\inf {\bf R}\Gamma_{\fm}(X,Y)=\sup {\bf R}\Hom_R({\bf R}\Hom_R(D,Y),X)$$ and $$-\sup {\bf R}\Gamma_{\fm}(X,Y)=\inf {\bf
R}\Hom_R({\bf R}\Hom_R(D,Y),X).$$
\end{enumerate}
\end{corollary}

The first part of the following corollary extends Suzuki's Duality Theorem and its second part extends
the Herzog-Zamani Duality Theorem.

\begin{corollary} Let $(R,\fm)$ be a Cohen Macaulay local ring possessing a canonical module $\omega$.
Let $M,N$  be two finitely generated $R$-modules and $i$ an integer.
\begin{enumerate}
\item[i)] Assume that one of the following conditions holds:\\
a) projective dimension of $M$ is finite,\\
b) both G-dimension of $M$ and injective dimension of $N$ are finite.\\
Then $$H_{\fm}^i(M,N)\cong \Hom_R(\Ext_R^{\dim R-i}(N,\omega\otimes_RM),E(R/\fm)).$$
\item[ii)] Assume that one of the following conditions holds:\\
a) injective dimension of $N$ is finite,\\
b) both Gorenstein injective dimension of $N$ and projective dimension of $M$ are finite. \\
Then $$H_{\fm}^i(M,N)\cong \Hom_R(\Ext_R^{\dim R-i} (\Hom_R(\omega,N),M),E(R/\fm)).$$
\end{enumerate}
\end{corollary}

\begin{prf} i) If G-dimension of $M$ is finite, then by \cite[Theorem 3.4.6]{C}, one has
$\omega\otimes_R^{{\bf L}}M\simeq \omega\otimes_RM$. Also, if Gorenstein injective dimension
of $N$ is finite, then \cite[Theorem 3.4.9]{C} asserts that ${\bf R}\Hom_R(\omega,N)\simeq
\Hom_R(\omega,N)$. Hence the conclusion is immediate by Theorem 3.3. Note that $\Sigma^{\dim R}
\omega$ is a normalized dualizing complex of $R$.
\end{prf}

\begin{example} None of Questions 1.1 and 1.2 have positive answers. To see this, let
$(R,\fm,k)$ be a non-regular Gorenstein local ring. Let $d:=\dim R$, and as before,
let $(-)^{\vee}$ denote the Matlis duality functor. Since $R$ is Gorenstein, by
\cite[Theorems 1.4.9 and 6.2.7]{C}, both G-dimension  of $k$ and Gorenstein injective
dimension of $k$ are finite. Assume that one of these questions has an affirmative answer.
Then, by Theorem 3.3, it turns out that
$$\Ext_R^i(k,k)\cong \underset{n}{\varinjlim}\Ext^{i}_{R}(k/\fm^{n}k,k)\cong  H^i_{\fm}(k,k)\cong \Ext^{d-i}_R(k,k)^{\vee}$$
for all non-negative integers $i$. This yields that $\Ext_R^i(k,k)=0$ for all $i\notin
\{0,1,\dots d \}$. So, $R$ is regular and we get a contradiction.
\end{example}

\section{Applications}

We start this section by proving a couple of lemmas.

\begin{lemma} Let $X\in \mathcal{D}_{\Box}^f(R)$ and $N\in \mathcal{D}_0^f(R)$.
Then  $$-\sup {\bf R}\Hom_R(N,X)=\inf \{\depth_{R_{\fp}}X_{\fp}|\fp\in \Supp_RN
\}=\depth_R(\Ann_RH_0(N),X).$$
\end{lemma}

\begin{prf}  The assertion follows immediately by \cite[Proposition 3.4]{Fo2} and
\cite[Proposition 2.10]{FI}.
\end{prf}

G-dimension of a complex $X\in \mathcal{D}^f_{\sqsupset}(R)$, $\Gdim_RX$, is defined by
$$\begin{array}{llll}\Gdim_RX&:=\inf\{\sup\{l\in \mathbb{Z}|Q_l\neq 0\}|Q \
\ \text{is a bounded to the right complex of} \\
&\text{totally reflexive R-modules and} \   Q\simeq X \}.
\end{array}
$$

\begin{lemma} Let $(R,\fm)$ be a local ring possessing  a normalized dualizing complex $D$,
$X\in \mathcal{D}_{\Box}^f(R)$ and $M,N$ two nonzero finitely generated $R$-modules.
\begin{enumerate}
\item[i)] If $M$ has finite G-dimension and $\Supp_RM\cap \Assh_RN\neq\emptyset$, then
$$\sup {\bf R}\Hom_R(N,D\otimes_R^{{\bf L}}M)\geq \dim_R N.$$
\item[ii)] If $N$ has finite Gorenstein injective dimension, then
$$\sup {\bf R}\Hom_R({\bf R}\Hom_R(D,N),X)=\depth R-\depth_R(\Ann_RN,X).$$
\end{enumerate}
\end{lemma}

\begin{prf} i) Let $\fp$ be a prime ideal of $R$. From \cite[15.17 c)]{Fo1} and
\cite[A.8.5.3]{C}, one has $\inf D_{\fp}=\dim R/\fp+\depth
R_{\fp}$. As $M\in \mathcal{A}^f(R)$, by \cite[Observation 3.1.7]{C},
it follows that $M_{\fp}\in \mathcal{A}^f(R_{\fp})$. So, by applying
\cite[Lemma A.6.4]{C} and \cite[A.6.3.2]{C}, we can deduce that
$$\begin{array}{llll} \depth_{R_{\fp}}M_{\fp}&=\depth_{R_{\fp}}
({\bf R}\Hom_{R_{\fp}}(D_{\fp},D_{\fp}\otimes_{R_{\fp}}^{{\bf L}}M_{\fp}))\\
&=\depth_{R_{\fp}}(D_{\fp}\otimes_{R_{\fp}}^{{\bf L}}M_{\fp})+\dim
R/\fp+\depth R_{\fp}.
\end{array}$$
Thus by Lemma 4.1 and the Auslander-Buchsbaum formula for G-dimension (see e.g.
\cite[Theorem 1.4.8]{C}), one has:
$$\begin{array}{llll} \sup {\bf R}\Hom_R(N,D\otimes_R^{{\bf L}}M)&
=-\inf\{-\dim R/\fp-\Gdim_{R_{\fp}}M_{\fp}|\fp\in \Supp_RN \}\\
&=\sup \{\dim R/\fp+\Gdim_{R_{\fp}}M_{\fp}|\fp\in \Supp_RN  \}\\
&=\sup \{\dim R/\fp+\Gdim_{R_{\fp}}M_{\fp}|\fp\in \Supp_RM\cap \Supp_RN  \}\\
&\geq \dim_RN.
\end{array}$$

ii) As $N\in \mathcal{B}^f(R)$, one has
$$N\simeq D\otimes_R^{{\bf L}}{\bf R}\Hom_R(D,N). \   \  (*)$$
By \cite[A.8.5.3]{C}, we have $\inf D=\depth R$. Set $s:=\depth
R$. Then applying Nakayama's Lemma for complexes (see e.g.
\cite[Corollary A.4.16]{C}) to $(*)$ yields that
$$\inf {\bf R}\Hom_R(D,N)=-\inf D=-s.$$ On the other hand, by
\cite[Proposition A.4.6]{C}, we have $$\sup {\bf R}\Hom_R(D,N)\leq
\sup N-\inf D=-s.$$ Hence $\Sigma^s{\bf R}\Hom_R(D,N)\in
\mathcal{D}_0^f(R)$. From $(*)$, one can conclude that $\Sigma^s{\bf
R}\Hom_R(D,N)$ and $N$ have the same support, and so Lemma 4.1
implies that
$$\begin{array}{llll} \sup {\bf R}\Hom_R({\bf R}\Hom_R(D,N),X)&=\sup
\Sigma^s{\bf R}\Hom_R(\Sigma^s{\bf R}\Hom_R(D,N),X)\\
&=s+\sup {\bf R}\Hom_R(\Sigma^s{\bf R}\Hom_R(D,N),X)\\
&=s-\depth_R(\Ann_RH_0(\Sigma^s{\bf R}\Hom_R(D,N)),X)\\
&=s-\depth_R(\Ann_RN,X).
\end{array}$$
\end{prf}

In the sequel, we establish a characterization of Cohen-Macaulay modules. It partially improves
\cite[Theorem 3.3]{HZ}. To this end, for a complex $Y\in \mathcal{D}_{\sqsubset}(R)$, we fix the
notation $Y^{\bot}$ for the full subcategory of $\mathcal{D}_0^f(R)$ whose objects are exactly those
complexes $X\in \mathcal{D}_{0}^f(R)$ for which $H_{\fm}^i(X,Y)=0$ for all $i\neq \depth_RY$.

\begin{proposition} Let $(R,\fm)$ be a local ring and $N$ a nonzero finitely generated $R$-module.
Consider the following conditions.
\begin{enumerate}
\item[i)] $N$ is Cohen-Macaulay.
\item[ii)] There is a nonzero $R$-module $M\in  N^{\bot}$ of finite
projective dimension such that $\Supp_RM\cap \Assh_RN\neq\emptyset$.
\item[iii)] There is a nonzero $R$-module $M\in  N^{\bot}$ of finite G-dimension such that $\Supp_RM\cap
\Assh_RN\neq\emptyset$.
\end{enumerate}
Then $i)$ and $ii)$ are equivalent and clearly $ii)$ implies $iii)$. In addition, if either projective
or injective dimension of $N$ is finite, then  all these conditions are equivalent.
\end{proposition}

\begin{prf} If a finitely generated $R$-module $M$ has finite G-dimension, then it is easy to check
that the $\widehat{R}$-module $\widehat{M}$ has finite G-dimension too. Also, if a finitely generated
$R$-module $M$ satisfies $\Supp_RM\cap \Assh_RN\neq\emptyset$, then $\Supp_{\widehat{R}}\widehat{M}\cap
\Assh_{\widehat{R}}\widehat{N}\neq\emptyset$. So, without loss of generality, we may and do assume that
$R$ is complete. So, $R$ possesses a normalized dualizing complex $D$.

$i)\Rightarrow ii)$ Assume that $N$ is Cohen-Macaulay. Then $H_{\fm}^i(R,N)=H_{\fm}^i(N)=0$ for all
$i\neq \depth_RN$, and so $R\in N^{\bot}$.

$ii)\Rightarrow iii)$ is clear.

Assume that either projective or injective dimension of $N$ is finite. We show $iii)$ implies $i)$.
Suppose that there exists a nonzero $R$-module $M\in N^{\bot}$ which has finite G-dimension.
Then $H_{\fm}^i(M,N)=0$  for all $i\neq \depth_RN$. Hence from Corollary 3.4 i) and Lemma 4.2 i),
we deduce that $$\depth_RN=-\inf {\bf R}\Gamma_{\fm}(M,N)=\sup {\bf R} \Hom_R(N,D\otimes_R^{{\bf L}}M)
\geq \dim_R N,$$ and so $N$ is Cohen-Macaulay.

$ii)\Rightarrow i)$ is similar to the proof of $iii)\Rightarrow i)$.
\end{prf}

Next, we establish the Gorenstein analogue of Proposition 4.3. It is worth to point out that it improves \cite[Proposition 3.5 1)]{HZ}. Recall that a non-homologically trivial complex $Y\in \mathcal{D}_{\Box}^f(R)$
is said to be Gorenstein if $\id_RY=\depth_RY$.

\begin{proposition} Let $(R,\fm,k)$ be a local ring and $Y\in \mathcal{D}_{\Box}^f(R)$ a non-homologically
trivial complex. The following are equivalent:
\begin{enumerate}
\item[i)] $Y$ is Gorenstein.
\item[ii)] $Y^{\bot}=\mathcal{D}_{0}^f(R)$.
\item[iii)] $k\in Y^{\bot}$.
\end{enumerate}
\end{proposition}

\begin{prf} $i)\Rightarrow ii)$  Let $X\in \mathcal{D}_{0}^f(R)$. By \cite[Theorem 2.7]{Y},
$\inf \{i\in \mathbb{Z}| H_{\fm}^i(X,Y)\neq 0\}=\depth_R Y$. Since $Y$ is Gorenstein, one has
$\id_RY=\depth_RY$, and so by Lemma 3.2, it turns out that $H_{\fm}^i(X,Y)=0$ for all $i\neq \depth_RY$.

$ii)\Rightarrow iii)$ is clear.

$iii)\Rightarrow i)$ Since $\Supp_RY\cap \Supp_Rk=\{\fm\}$, by \cite[Lemma 2.4]{Y}, one has $H_{\fm}^i(k,Y) =\Ext_R^i(k,Y)$ for all integers $i$. Thus $\Ext_R^i(k,Y)=0$ for all $i\neq \depth_RY$. By \cite[A.5.7.4]{C},
this yields that $\id_RY=\depth_RY$.
\end{prf}

\begin{lemma} Let $(R,\fm)$ be a local ring possessing a normalized dualizing complex $D$ and $X\in \mathcal{D}_{\Box}^f(R)$. Assume that $X$ has finite G-dimension. Then
$$\depth_R(\fa,X)-\dim R\leq \depth_R(\fa,D\otimes_R^{{\bf L}}X)\leq \depth_R(\fa,X)-\depth R$$
for all ideals $\fa$ of $R$.
\end{lemma}

\begin{prf} Let $\underline{\fa}$ be a generating set for a given ideal $\fa$ of $R$. As
$\Check{C}(\underline{\fa})$ is a bounded complex of flat $R$-modules, $X\in \mathcal{A}^f(R)$ and
$\inf D=\depth R$,  \cite[Theorem 4.7 i)]{CH} and \cite[Proposition 3.3.7 a)]{C}
imply that
$$\sup {\bf R}\Gamma_{\fa}(X)+\depth R\leq \sup {\bf R}\Gamma_{\fa}
(D\otimes_R^{{\bf L}}X)=\sup (D\otimes_R^{{\bf L}}{\bf R}\Gamma_{\fa}(X))\leq \sup {\bf
R}\Gamma_{\fa}(X)+\dim R.$$
But for any complex $Z\in \mathcal{D}_{\Box}^f(R)$, by \cite[Proposition 3.14 c)]{Fo2}, one has
$\depth_R(\fa,Z)=-\sup {\bf R}\Gamma_{\fa}(Z)$. This completes the proof.
\end{prf}

The following result can be considered as Grothendieck's non-vanishing Theorem in the context of generalized
local cohomology modules. It also improves \cite[Theorem 3.5]{DH}.

\begin{proposition} Let $(R,\fm)$ be a local ring, $N$ a nonzero finitely generated $R$-module and $X\in
\mathcal{D}_{\Box}^f(R)$.
\begin{enumerate}
\item[i)] Assume that one of the following conditions is satisfied:\\
a) projective dimension of $X$ is finite,\\
b) G-dimension of $X$ and projective dimension of $N$ are finite.\\
Then $$\depth R-\depth(\Ann_RN,X)\leq \cd_{\fm}(X,N)\leq \dim R-\depth(\Ann_RN,X).$$
\item[ii)] Assume that one of the following conditions is satisfied:\\
a) injective dimension of $N$ is finite, \\
b) Gorenstein injective dimension of $N$ and either projective dimension or injective dimension of $X$ are finite.\\
Then $$\cd_{\fm}(X,N)=\depth R-\depth(\Ann_RN,X).$$
\end{enumerate}
\end{proposition}

\begin{prf} Let $Z\in \mathcal{D}_{\Box}^f(R)$. Clearly, then one has $Z\otimes_R\widehat{R}\in \mathcal{D}_{\Box}^f(\widehat{R})$. If projective (resp. injective) dimension of $Z$ is finite, then
$Z\otimes_R\widehat{R}$ has finite projective (resp. injective) dimension over $\widehat{R}$. Also, it
is easy to check that if G-dimension of $Z$ is finite, then so is G-dimension of $Z\otimes_R\widehat{R}$
over $\widehat{R}$. By \cite[Theorem 3.6]{FF} if Gorenstein injective dimension of $Z$ is finite, then
so is Gorenstein injective dimension of $Z\otimes_R\widehat{R}$ over $\widehat{R}$. On the other hand, since $\widehat{R}$ is a faithfully flat $R$-module, for any complex $W$,  one has
$\sup W=\sup (W\otimes_R\widehat{R})$ and $\inf W=\inf (W\otimes_R\widehat{R})$. Hence
$$\begin{array}{llll} \depth_{\widehat{R}}(\Ann_{\widehat{R}}\widehat{N},
X\otimes_R\widehat{R})&=\depth_{\widehat{R}}((\Ann_RN)\widehat{R},X\otimes_R
\widehat{R})\\
&=-\sup {\bf R}\Hom_{\widehat{R}}(\widehat{R}/(\Ann_RN)\widehat{R},X\otimes_R
\widehat{R})\\
&=-\sup ({\bf R}\Hom_R(R/\Ann_RN,X)\otimes_R\widehat{R})\\
&=-\sup {\bf R}\Hom_R(R/\Ann_RN,X)\\
&=\depth(\Ann_RN,X).
\end{array}
$$
Similarly, one has $\cd_{\fm \widehat{R}}(X\otimes_R\widehat{R},\widehat{N})=\cd_{\fm}(X,N)$. Thus,
we may and do assume that $R$ is complete. Hence $R$ possesses a normalized dualizing
complex $D$. In case i), the G-dimension of $X$ is finite. By using Corollary 3.4 i) and Lemma 4.1,
we can deduce that:
$$\cd_{\fm}(X,N)=\sup {\bf R}\Hom_R(N,D\otimes_R^{{\bf L}}X)
=-\depth_R(\Ann_RN,D\otimes_R^{{\bf L}}X).$$ Hence Lemma 4.5 completes the proof of i).

In case ii), the Gorenstein injective dimension of $N$ is finite. Hence, the conclusion follows by Corollary
3.4 ii) and Lemma 4.2 ii).
\end{prf}

The following result partially generalizes the Intersection inequality.

\begin{proposition} Let $(R,\fm)$ be a local ring and $M,N$ two nonzero finitely generated $R$-modules
such that $ \Supp_RM\cap \Assh_RN\neq\emptyset$. Assume that one of the following conditions is satisfied:
\begin{enumerate}
\item[i)] projective dimension of $M$ is finite.
\item[ii)] both G-dimension of $M$ and projective dimension of $N$ are finite.
\item[iii)] injective dimension of $N$ is finite.
\item[iv)] both Gorenstein injective dimension of $N$ and injective dimension of $M$ are finite.
\end{enumerate}
Then $$\dim_R N\leq \dim_R{\bf R}\Hom_R(M,N)\leq -\inf {\bf R}\Hom_R(M,N)+\dim_R(M\otimes_RN).$$
\end{proposition}

\begin{prf} It is easy to check that $\dim_R{\bf R}\Hom_R(M,N)=\dim_{\widehat{R}}{\bf
R}\Hom_{\widehat{R}}(\widehat{M},\widehat{N})$. Hence, in view of the proof of Proposition 4.6, we may
assume that $R$ is complete. So, $R$ possesses a normalized dualizing complex. In each case, it
follows that ${\bf R}\Hom_R(M,N)\in \mathcal{D}_{\Box}^f(R)$. By Grothendieck's non-vanishing Theorem
\cite[Proposition 3.14 d)]{Fo2}, one has $\cd_{\fm}(M,N)=\dim_R{\bf R}\Hom_R(M,N)$. Thus \cite[Corollary
3.2]{DH} yields the right hand inequality. In cases i) and ii), the left hand inequality follows by
Corollary 3.4 i) and Lemma 4.2 i). Let $\fp_0\in  \Supp_RM\cap \Assh_RN\neq\emptyset$. In each of the
cases iii) and iv), our assumptions yield that $R$ is Cohen-Macaulay, and so one has:
$$\begin{array}{llll} \dim_RN&=\dim R/\fp_0\\
&=\dim R-\Ht \fp_0\\
&\leq \dim R-\depth_{R_{\fp_0}}M_{\fp_0}\\
&\leq \dim R-\inf\{\depth_{R_{\fp}} M_{\fp}|\fp\in \V(\Ann_RN)\}\\
&=\dim R-\depth_R(\Ann_RN,M).
\end{array}
$$
Hence in these cases,  the left hand inequality follows by Proposition 4.6 ii).
\end{prf}


\end{document}